%% ok, attempting to put this on arXiv, february 2026 

% First version written 15-jul-1988
% Intention, as of August 1989: Send it off to Statistics
% & Probability Letters, with others, from Nils 1976 Tromsoe.
% now it's 2021, tenk det; and july 2022, tenk det ogsaa  
% topoffile

%% sent to Stat and Prob Letters july 2022,
%% https://www.editorialmanager.com/stapro/default1.aspx
%% where i need to log in as ... with slavatrudu444 

\magnification\magstep1
\baselineskip=13pt
\vsize=24.0truecm

\font\csc=cmcsc10
\font\smallrm=cmr9
\font\smallsl=cmsl9
\font\smallbf=cmbx9

\def\d{{\rm d}}
\def\E{{\rm E}}

\centerline{\bf A simple proof of the discreteness of Dirichlet processes}

\bigskip 
\centerline{\bf Nils Lid Hjort}

\medskip
%% \centerline{\sl Norwegian Computing Centre {\it\&} University of Oslo}
\centerline{\sl Department of Mathematics, University of Oslo}

\medskip
\centerline{\sl -- February 2026 --}
%% \centerline{\sl -- December 1976, but now it's July 2022 --}

{{\medskip\narrower\narrower
\noindent {\csc Abstract.} 
That Dirichlet processes are discrete with probability 1 
is demonstrated once more. And yes, these two pages
spent fifty years in Norwegian.

\smallskip
\noindent {\csc Key Words:} Dirichlet processes, discreteness
\smallskip}} 

\medskip\noindent
Let $\cal X$ be some sample space, 
and let $P$ be a Dirichlet process with parameter $kP_0$ on $\cal X$. 
Here $k$ is a positive scalar and $P_0$ is a probability
measure on $\cal X$. Thus $P$ is a random element of $\cal M$, the set
of all probability measures on the sample space. There are
various proofs in the literature demonstrating that such a 
$P$ is discrete with probability one; 
a possibly partial list includes Ferguson (1973, 1974),
Blackwell (1973), Blackwell and MacQueen (1973), 
Doksum (1974), Kingman (1975), 
Berk and Savage (1979), 
Hjort (1986), Krasker and Pratt (1986), 
Sethuraman and Tiwari (1982), and Sethuraman (1994). 
In particular, the so-called stick-breaking representation
of a Dirichlet process, due to the two latter references,
spell out in a constructive fashion that the random $P$
must indeed be discrete. 
Below follows yet another proof, from Hjort's 1976 thesis (p.~18).

For a general overview of contemporary nonparametric
Bayesian statistics, which includes the frequent use
of Dirichlet processes in various forms, see
Hjort (2003), and several chapters of the 
Hjort, Holmes, M\"uller, Walker (2010) book. 

A useful lemma concerning the expected value of a function $f({P},X)$,
where $X$ is sampled from the random $P$, states that
$$\E\,f({P},X)
        =\int_{\cal M}\int_{\cal X}f(P,x)\,P(\d x)\,{\cal D}_{kP_0}(\d P)
        =\int_{\cal X}\int_{\cal M}f(P,x)\,
        {\cal D}_{kP_0+\delta_x}(\d P)\,P_0(\d x).$$
Some measure theoretic details must be tended to here: there is some
sigma-field $\cal A$ on $\cal X$; $\cal M$ is equipped, for example,
by the Borel sets determined by set-wise convergence; and $f$ must be
measurable in $(P,x)$. Also, ${\cal D}$, appropriately subscripted, 
is used to denote $P$'s probability distribution, 
and $\delta_x$ denotes the probability measure with unit mass 
at position~$x$. 

The lemma is related in an obvious way to two well-known 
facts about the Dirichlet process,
namely that $X$ as above has unconditional distribution $P_0$, 
and that $P$, conditionally on an observed $X=x$, is 
Dirichlet with updated parameter $kP_0+\delta_x$. 
There is a natural extension 
to functions $f({P},X_1,\ldots,X_n)$. 
The lemma was proved and used for various causes in Hjort (1976), 
and has later on been re-discovered on appropriate occasions; 
see Lo (1984) for but one example. 

Introduce $A_P=\{x\colon P\{x\}>0\}$, the set of atoms for a given $P$, 
and define
$$\eqalign{
  H_\gamma(P)
  &=\E_P P\{X\}^\gamma=\int P\{x\}^\gamma\,\d P(x)
   =\sum_{x\in A_P} P\{x\}^{\gamma+1}, \cr
H_0(P)&=\lim_{\gamma\rightarrow 0^+}H_\gamma(P)=\sum_{x\in A_P}P\{x\}. \cr}$$
A $P$ is discrete if and only if $H_0(P)=1$. Employ the lemma to get
$$\eqalign{ \E\,H_\gamma({P})=\E\,{P}\{X\}^\gamma
        &=\int_{\cal X}
        {\Gamma(kP_0\{x\}+1+\gamma)\over \Gamma(kP_0\{x\}+1)}
        {\Gamma(k+1)\over \Gamma(k+1+\gamma)}\,P_0(\d x)\cr
        &\ge{\Gamma(1+\gamma)\over \Gamma(1)}
        {\Gamma(k+1)\over \Gamma(k+1+\gamma)}.}$$
From this and $0\le H_\gamma(P)\le H_0(P)\le 1$ follows 
$\E\,H_0({P})=1$ and 
$H_0({P})=1$ with ${\cal D}_{kP_0}$-probability one. The single
measure theoretic caveat here is that $({\cal X},{\cal A})$ must be such
that $P\{x\}$ is simultaneously measurable. It suffices that $\cal A$ 
is the set of Borel sets from a metric which makes the sample 
space separable. Such conditions make $H_\gamma$ measurable in $P$,
and also entails the measurability of the set ${\cal M}_0$ of all
discrete probability measures. 

Arguments similar to those above show that each set $A$ with positive
$P_0$-measure must have positive ${P}$-atoms, with probability one.  
If in particular $P$ is Dirichlet with parameter $kP_0$ 
on the real line, with random distribution function $F$, then $F$
has infinitely many jumps on each interval with positive $P_0$-measure.

\bigskip
{\bf Acknowledgments.} 
The work reported on in this note has spent fifty years
in Norwegian, and is based on pages 18--19 of my graduate thesis,
written while a lecturer 
%% an amanuensis 
in Troms\o{} 1976.
During that time I benefited from discussions with
Berit Sandstad and Tore Schweder.

\bigskip
\centerline{\bf References}

% \def\ref#1{{\noindent\hangafter=1\hangindent=20pt
%  #1\smallskip}}          
% \parindent0pt
% \baselineskip11pt
% \parskip3pt 
% \medskip 

\parindent0pt\smallrm
\medskip
Basu, D.~and Tiwari, R.C. (1982).
A note on the Dirichlet process.
In {\smallsl Statistics and Probability: Essays in Honour of C.R.~Rao},
eds.~Kallianpur, Krishnaiah, Ghosh, 89--103.

Berk, R.H.~and Savage, I.R. (1979). 
Dirichlet processes produce discrete measures: an elementary proof.
{\smallsl Contributions to Statistics: Jaroslav H\'ajek Memorial Volume},
Academia, North-Holland, 25--31.

Blackwell, D. (1973).
Discreteness of Ferguson selections. 
{\smallsl Annals of Statististics} {\smallbf 1}, 356--158. 

Blackwell, D.~and MacQueen, J.B. (1973).
Ferguson distributions via P\'olya urn schemes.
{\smallsl Annals of Statistics} {\smallbf 1}, 353--355. 

Doksum, K. (1974). 
Tailfree and neutral random probability distributions and their
posterior distributions.
{\smallsl Annals of Probability} {\smallbf 2}, 183--201.

Ferguson, T.S. (1973).
A Bayesian analysis of some nonparametric problems. 
{\smallsl Annals of Statistics}~{\smallbf 1}, 209--230.

Ferguson, T.S. (1974).
Prior distributions on spaces of probability measures.
{\smallsl Annals of Statistics}~{\smallbf 2}, 615--629. 

Hjort, N.L. (1976). 
{\smallsl Applications of the Dirichlet process to some nonparametric
estimation problems} (in Norwegian), graduate thesis, 
University of Troms\o. 
Abstracted in {\smallsl Scandinavian Journal of Statistics} 1977.

Hjort, N.L. (1986). 
Contribution to the discussion of Diaconis and 
Freedman's `On the consistency of Bayes estimates'. 
{\smallsl Annals of Statistics} {\smallbf 14}, 49--55.

Hjort, N.L. (2003).
Topics in nonparametric Bayesian statistics
(with discussion). 
In {\smallsl Highly Structured Stochastic Systems}
(eds.~P.J.~Green, N.L.~Hjort and S.~Richardson),
Oxford University Press. 

Hjort, N.L., Holmes, C.C., M\"uller, P., and Walker, S.G. (2010).
{\smallsl Bayesian Nonparametrics}. 
Cambridge University Press, Cambridge. 

Kingman, J.F.C. (1975). 
Random discrete distributions (with discussion). 
{\smallsl Journal of the Royal Statistical Society Series B} 
{\smallbf 37}, 1--22. 

Krasker, W.S.~and Pratt, J.W. (1986). 
Contribution to the discussion of Diaconis and 
Freedman's `On the consistency of Bayes estimates'. 
{\smallsl Annals of Statistics} {\smallbf 14}, 42--45.

Lo, A.Y. (1984). 
On a class of Bayesian nonparametric estimates: I. Density estimates.
{\smallsl Annals of Statistics} {\smallbf 12}, 351--357. 

Sethuraman, J.~and Tiwari, R. (1982).
Convergence of Dirichlet measures and the interpretation of their
parameter. In {\smallsl Proceedings of the Third
Purdue Symposium on Statistical Decision Theory and Related Topics}
(eds.~S.S.~Gupta and J.~Berger), 305--315. Academic Press, New York.

Sethuraman, J. (1994).
A constructive definition of Dirichlet priors. 
{\smallsl Statistica Sinica} {\smallbf 4}, 639--650. 

\bye